\documentclass[11pt,letterpaper]{amsart}

\usepackage[utf8]{inputenc}
\usepackage{amssymb, amscd, amsmath, epsfig, mathtools}
\usepackage{amsthm}
\usepackage{enumerate}
\usepackage{xcolor}
\usepackage{scalerel}
\usepackage{soul}
\usepackage{tikz-cd}
\usepackage{esint}
\usepackage{cancel}
\usepackage{hyperref}
\usepackage{slashed}
\usepackage{url}
\usepackage[normalem]{ulem}

\allowdisplaybreaks

\usepackage[margin=1.0in]{geometry}

\newtheorem{theorem}{Theorem}[section]
\newtheorem*{theorem*}{Theorem}

\newtheorem{corollary}{Corollary}[section]
\newtheorem*{corollary*}{Corollary}
\newtheorem{lemma}{Lemma}[section]
\newtheorem{proposition}{Proposition}[section]
\theoremstyle{definition}
\newtheorem{remark}{Remark}[section]

\newenvironment{introtheorem}[1]{%
  \renewcommand\thetheorem{#1}
  \theorem
}{\endtheorem}

\newcommand{\C}{\mathbb C}
\newcommand{\R}{\mathbb R}

\newcommand{\soo}{{\mathbb S}^1}
\newcommand{\calC}{\mathcal C}
\newcommand{\calW}{\mathcal W}
\newcommand{\calL}{\mathcal L}
\newcommand{\dvol}{ \text{dvol}_{g}}

\newcommand{\dvols}{\text{dvol}_{\sigma^*g}}

\begin{document}

\title[$ \mathbb{S}^{1} $-Stability Conjecture]{A Codimension  Two  Approach to The $ \mathbb{S}^{1} $-Stability Conjecture}
\author[S. Rosenberg and J. Xu]{Steven Rosenberg and Jie Xu}
\address{
Department of Mathematics and Statistics, Boston University, Boston, MA, U.S.A.}
\email{sr@math.bu.edu}
\address{
Department of Mathematics, Northeastern University, Boston, MA, U.S.A.}
\email{jie.xu@northeastern.edu}

\date{}							

\begin{abstract} J.~Rosenberg's $\soo$-stability conjecture states that a closed oriented manifold $X$ admits a positive scalar curvature metric iff $X\times \soo$ admits a positive scalar curvature metric $h$.  As pointed out by J.~Rosenberg and others, there are known counterexamples in dimension four. Using conformal geometry techniques, we prove this conjecture whenever 
$h$ satisfies a geometric bound which measures the discrepancy between $\partial_\theta\in T\soo$ and the normal vector field to $X\times \{P\}$, for a fixed $P\in \mathbb{S}^1.$
\end{abstract}

\maketitle

\section{Introduction}
In \cite[Conj.~1.24]{JR}, J.~Rosenberg proposed the
$ \mathbb{S}^{1} $-stability conjecture for metrics of 
positive scalar curvature (PSC).
\medskip

\noindent {\bf{$\soo$-Stability Conjecture}}. Let $ X $ be a closed oriented manifold. Then $ X $ admits a PSC metric if and only if $ X \times \mathbb{S}^{1} $ admits a PSC metric.
\medskip

In \cite[Rmk.~1.25]{JR},  Rosenberg gave a counterexample in dimension four of an odd degree hypersurface in $\C\mathbb{P}^3$, based on  Seiberg-Witten theory. Other counterexamples in dimension four are in \cite[Rmk.~5]{KS}. The conjecture has been verified in dimensions  $2,  3, 5, 6 $  and
for a large class of spin manifolds \cite[Cor.~1.5]{Zeidler}, using
a combination of minimal surface methods and band-width estimates.   (See {\em e.g.,}  
\cite[\S1.1]{carlotto} for the ingredients for dimension 2, 
\cite[Cor.~7.34]{Chodosh} for dimension 3, and \cite[Rmk.~2.26]{Rade} for dimensions 5, 6.)

The nontrivial direction of the conjecture states that if $(X\times \soo, h)$ is a PSC metric, then $X$ admits a PSC metric. 
We introduce a PDE approach to solve the $ \mathbb{S}^{1} $-stability conjecture in all dimensions under a  restriction on $h$.  To state the main theorem, 
let $\partial_\theta$ be the usual tangent vector field to 
$\soo$, and
let $\mu$ be the unit normal vector field to the slice $X\times \{P\}$ for a fixed $P\in \soo$,
where $\mu$ is chosen so that $h(\mu,\partial_\theta)>0.$
Define the $h$-angle $\angle_h(\mu, \partial_\theta)$ by $h(\mu,\partial_\theta) = \Vert\mu\Vert_h\Vert\partial_\theta\Vert_h \cos(\angle_h(\mu, \partial_\theta)) = \Vert\partial_\theta\Vert_h \cos(\angle_h(\mu, \partial_\theta))$.

\medskip

\begin{introtheorem}{3.1}\label{intro:thm}
Let $ X $ be an oriented closed manifold with $ \dim X  \geqslant 2 $. If $ X \times \mathbb{S}^{1} $ admits a PSC metric $ h $ such that for some $P\in \soo$, we have
\begin{equation}\label{eq:intro}
\angle_h(\mu, \partial_\theta) < \frac{\pi}{4} 
\end{equation}
on $X\times \{P\},$
then $ X \approx X \times \{P\}$ admits a PSC metric.
\end{introtheorem} 
To be precise, the $\soo$-conjecture holds on $X\times \{P\}$ if $\partial_{\mu}$ lies inside the cone with axis $\partial_{\theta}$ and $h$-angle $\pi/4$ on each $T_{(x,P)}(X\times\soo)$, for all $x\in X.$

\medskip

The main techniques used are a mix of the Gauss-Codazzi equation and conformal geometry,
along with the introduction of an extra dimension.
To outline the proof of the main theorem, we refer to the diagram (\ref{diagram}). We take the
PSC metric $h$ on $X\times \soo$ and extend it to a product metric $g$ on $M := W\times \soo := X\times [-1,1]\times \soo$. Note that  $X$ has codimension two in $M$.  This induces a metric $\sigma^*g$ on $W$, which we use to construct a PDE (\ref{CG:eqn1}), which is elliptic if (\ref{eq:intro}) holds.
This PDE is related to a conformal transformation of the Gauss-Codazzi equation, and the choice of the inhomogeneous term $F$ 
in (\ref{CG:eqn1}) is crucial. Using the ellipticity, we find a 
smooth solution $u$ to  
(\ref{CG:eqn1}) on $W$, such that $u_W := u+1$ is positive. 
We pull back $u_W$ to $u_M$ on $M$, and  use $u_M$ to define a conformal transformation of $g$ to $\tilde g = u_M^{\frac{4}{n-2}}g$.  Restricting $\tilde g$ to $X\times S^1$ and then to $X$ gives a PSC metric $\tau^*\imath^*\tilde g$ on $X$.  

It is natural to ask why the extra dimension $[-1,1]$ is introduced.  If we take a conformal transformation $\tilde h = u^{\frac{4}{n-2}}h$ on $X\times S^1$ and use the Gauss-Codazzi equation to try to prove that the scalar curvature of $X\times \{P\}$ is positive, we find that $u$ must satisfy a non-elliptic PDE.  In particular, we cannot guarantee the existence of a smooth solution with the necessary small $\calC^{1,\alpha}$ estimate (see Prop.~\ref{CG:lemma2}). 
By adding in the extra dimension, we produce and solve 
(\ref{CG:eqn1}), provided (\ref{eq:intro}) holds.  This is discussed before Lem.~\ref{CG:lemma6}.

As an outline of the paper, in \S2 we prove a series of technical results.  
The  
elliptic PDE (\ref{CG:eqn1}) on $W$ related to the conformal transformation $g\mapsto \tilde g$ is given in 
Prop.~\ref{CG:lemma2}. We prove that (\ref{CG:eqn1}) has a smooth solution $u$ with  $\Vert u\Vert_{\calC^{1,\alpha}}$ arbitrarily small,   by choosing the inhomogeneous term $F$  
to be concentrated near $X = X\times \{0\}.$  
In particular, Lem.~\ref{CG:second} gives a partial $ \calC^{2} $-estimate for $u_W$, where the 
hypothesis that $h$ is PSC is essentially used. We also prove several technical lemmas which allow us to compare the Laplacians and gradients on $M, W,$ and $X\times \soo$.

In \S3 we prove the main theorem Thm.~ \ref{intro:thm}, essentially by using Gauss-Codazzi on $\tilde g$ twice to move from $M$ back to $X\times \soo$ and then back to $X$. The technical results from \S2 appear in the Gauss-Codazzi formulas and various conformal transformation rules in the proof of Thm.~ \ref{intro:thm}.

We thank B.~Hanke, D.~Ruberman, T.~Schick and B.~Sen for corrections to an earlier version of this paper.

\section{Technical Preliminaries} 
In this section, we prove preliminary results that are used in the proof of Thm.~
\ref{MAIN:thm1} in \S3.

We introduce the setup and notation. 
$X$ is an oriented closed manifold such that $ X \times \mathbb{S}^{1} $ has a PSC metric $ h $.  We assume dim$(M) = n-1\geq 2.$
Set 
$$Y = X\times \soo,\ W = X \times [-1, 1],\  M =  W \times \mathbb{S}^{1} = Y\times [-1, 1].$$ 
  $ X \times \mathbb{S}^{1} 
=X\times \{0\} \times \soo$ is an embedded hypersurface of $M$.

Put the product metric $ g = h \oplus dt^{2} $ on $ M $, for $ t \in [-1, 1]. $ 
For a fixed $P\in \mathbb{S}^1$, define $ \sigma : W \rightarrow W \times \mathbb{S}^{1}$ by
$\sigma(w) = (w, P),$ and $\tau: X \rightarrow X \times \mathbb{S}^{1}$, $\tau(x) = (x,P).$
For the inclusion maps  
$ \imath  = \imath_0: X\rightarrow W$  and 
$\imath = \imath_0: X \times \mathbb{S}^{1} \rightarrow M $,  where $\imath_0(x) = (x,0)$ and
$\imath_0(x,P) = (x,0,P)$,
we have induced metrics $ \sigma^{*} g $ on $ W $ and $ \tau^{*}\imath^{*} g = \tau^*h $ on $ X $.
The following diagram summarizes the setup.
\begin{equation}\label{diagram}
\begin{tikzcd} (Y = X\times \soo, h=\imath^*g) \arrow[r,"\imath"] & (M= W\times \soo = Y\times [-1,1],g) \\
(X,\tau^* \imath^*g= \imath^*\sigma^*g) \arrow[r,"\imath"]\arrow[u,"\tau"] &  (W = X \times [-1, 1], \sigma^*g) 
\arrow[u,"\sigma"]
\end{tikzcd}
\end{equation}
We will frequently identify $x\in X$ with $ (x,0)\in W$, $y:=(x,P)\in Y$, and $(x,0,P)\in M.$ Note that  
dim $M\geq 4.$  

Let $ \theta $ be the usual coordinate on $\soo$ with corresponding global vector field $\partial_\theta$ on $W$ and $M$. Let $\mu$ be the outward normal vector field to $X\times \{Q\}$, for all $Q\in \soo,$ chosen so that  $h(\mu,\partial_\theta) >0.$
$\mu$ is a global vector field on $W$.  Locally on $W$, we can consider $\mu$ to be the tangent vector $\partial_{x^n}$, where $(x^1,\ldots,x^{n-1})$ are local coordinates on $U\subset X$.  In these coordinates, $ h(\mu, \mu) = h_{nn} = g_{nn} : = g_{\mu \mu} $, where the expressions in $g$ are valid on $U\times [-1, 1].$
Below, we will drop some inclusion maps and  write $g_{nn}$ on  
both $ M $ and
$ X\times \lbrace 0 \rbrace  \times \mathbb{S}^{1}.$
\medskip

As 
outlined in the Introduction, in this section we  solve a  
PDE (\ref{CG:eqn1})
involving $\sigma^*g$  on $W$ (Prop.~\ref{CG:lemma2}). 
This PDE is elliptic under the condition (\ref{eq:intro}) (Lem.~ \ref{CG:lemma6}).
We modify a solution $u$ to a positive function $u_W$ on $W$, pull back $u_W$ to $u_M$ on $M$,
and use $u_M$ to conformally change $g$ to $\tilde g = u_M^{\frac{4}{n-2}}g$.  The main theorem in \S3 proves that for a careful choice of the inhomogeneous term $F$ in (\ref{CG:eqn1}), $\tau^*\imath^*\tilde g$ is PSC, provided 
(\ref{eq:intro}) holds.
To this end,
 the diagram indicates that we have to compare Laplacians $\Delta_g$, $ \Delta_{\imath^{*}g} $ and $\Delta_{\sigma^*g}$ (Lem.~\ref{CG:lemma3}), the scalar curvatures of $ \tilde{g}$ and $\imath^* \tilde{g}$,
 and the scalar curvatures of $\imath^{*} \tilde g$ and $\tau^{*} \imath^*\tilde g$ (see the proof of Thm.~ \ref{MAIN:thm1}) via the Gauss-Codazzi equation.

\medskip

The function $u_M$ defined above is independent of $\theta.$  However, comparing the scalar curvatures for induced metrics via Gauss-Codazzi requires the use of the normal vector $\mu.$  For the product metric on $X\times S^1$, $\mu = \partial_\theta$, so in the non-product case we have to compare these vector fields. On $X\times S^1$,  set
\begin{equation}\label{eq:Vdef} \mu = a\partial_\theta + \sum_{i=1}^{n-1} b^i\partial_i := a\partial_\theta + V,
\end{equation}
for a  local coordinate frame $\{\partial_i\}$ on $X$.  
Note that
$$ 1 = h(\mu,\mu) = h(a\partial_\theta,\mu) + b^ih(\partial_i, \mu) = a h(\partial_\theta,\mu) 
\Rightarrow a = h(\partial_\theta,\mu) ^{-1},$$
where clearly $h(\partial_\theta,\mu)\neq 0.$
Since $\mu, a\partial_\theta $ are globally defined vector fields, so is the vector field $V$
on $X \times \{P\}$ for each $P\in \soo.$

Finally, $\mu, \partial_\theta$  and hence $V$ pull back to  $t$-independent vector fields on $M$. Since $\sigma = \sigma_P$ identifies $X\times [-1, 1]\times \{P\}$ with $W$,  it is easy to check that $|V_{(x,P)}|^2_h
= |V_{(x,t)}|^2_{\sigma^*g}$ for all $t$.

\medskip

The operator $L'$ in the next lemma appears in the proof of the main theorem; see 
(\ref{eq:22}).  This lemma
motivates
the introduction of the extra dimension $[0,1]$ in $W$; if we work on $X\times \soo$, the  operator  analogous to $L'$ associated to a conformal change of $h$ would be 
$\nabla_{\mu}\nabla_{\mu} -  \Delta_{h}$, which is never elliptic.

Note that  (\ref{eq:htn}) is equivalent to the angle hypothesis (\ref{eq:intro}) in Thm.~ \ref{intro:thm}.

\begin{lemma}\label{CG:lemma6} Fix $P\in \soo,$ and identify $ X \times [-1, 1] \times \lbrace P \rbrace$ with $ W$. 
If 
\begin{equation}\label{eq:htn} \frac{h(\partial_\theta,\partial_\theta)}{h(\mu,\partial_\theta)^{2}}<
2,\end{equation}
 then 
the operator
\begin{equation*}
 L':=\nabla_{V}\nabla_{V} -  \Delta_{\sigma^{*}g}
\end{equation*}
is elliptic on $ W $.
\end{lemma}

We use the convention that $-\Delta_\ell$ is positive definite for a metric $\ell.$

\begin{proof}
In local coordinates,
 \begin{equation*}
    \nabla_{V}\nabla_{V} = \sum_{i, j = 1}^{n - 1} b^{i}(x) b^{j}(x) \frac{\partial^{2}}{\partial x^{i} \partial x^{j}} + G_{1}(x_{0})
\end{equation*}
where $ G_{1}(x) $ is a linear first order operator. In 
 Riemannian normal coordinates $(x^1,\ldots,x^{n-1})$ centered at a fixed $x_0\in X$, we have
\begin{equation*}
    \Delta_{\sigma^{*}g}|_{x_{0}} = \sum_{i = 1}^{n - 1} \frac{\partial^{2}}{\partial \left( x^{i} \right)^{2} } + \frac{\partial^{2}}{\partial t^{2}}.
\end{equation*}
The principal symbol of $L'$ at $ x_{0} $ is
\begin{equation*}
 \sigma_2(L')(x_0, \xi) =   - \sum_{i, j = 1}^{n - 1} b^{i}(x_{0})b^{j}(x_{0})
      \xi^{i} \xi^{j} 
     + \sum_{i = 1}^{n - 1} \left(\xi^{i} \right)^{2} + (\xi^n)^2.
\end{equation*}
Proving $ \sigma_2(L')(x_0, \xi)>0$ for $\xi\neq 0$ is equivalent to proving that
\begin{equation*}
    B : = - 
    \begin{pmatrix} (b^{1})^{2} & b^{1}b^{2} &\dotso & b^{1}b^{n - 1} \\ b^{2}b^{1} &( b^{2})^{2} & \dotso & b^{2}b^{n - 1} \\ \vdots & \vdots & \ddots & \vdots \\ b^{n - 1}b^{1} & b^{n - 1}b^{2} & \dotso & (b^{n - 1})^{2} \end{pmatrix} + I_{n - 1} : = B' + I_{n - 1}
\end{equation*}
is positive definite. 
Sylvester's criterion   states 
 the symmetric matrix $ B $ is positive definite iff all its $ k \times k $ principal minors 
 $ B_{k} $ have positive determinants.

 We claim that
\begin{equation}\label{CG:eqnR1}
 \det(B_{k}) = 1 + {\rm Tr}(B_{k}') = 1 -  \sum_{i = 1}^{k} (b^{i})^{2}, k = 1, \dotso, n - 1.
\end{equation}
We show this for $ B_{n-1} = B $, as the argument for the other $ B_{k} $ 
is identical. We first observe that if  $ b^{i} \neq 0, \forall i, $ then all rows of $ B' $ are proportional to each other. It follows that the kernel of $B'$ has dimension $ n - 2 $,  so $0$ is an eigenvalue of  $ B' $ of multiplicity $ n - 2 $. Since the sum of  the eigenvalues of $ B' $ is  the trace $- \sum_{i = 1}^{n - 1} (b^{i})^{2} $, 
the nontrivial eigenvalue must be $- \sum_{i = 1}^{n - 1} (b^{i})^{2} $.
Thus the eigenvalues of $B$ are 
\begin{equation*}
    \lambda_{1} = - \sum_{i = 1}^{n - 1} (b^{i})^{2} + 1, \lambda_2 = \ldots = \lambda_{n-1} = 1,
\end{equation*}
which proves (\ref{CG:eqnR1}). 
 If some  $ b^{i} $ vanish, then the claim reduces to a lower dimensional case by removing the corresponding  rows and columns of all zeros. 
 
Therefore, $L'$ is elliptic if 
 \begin{equation}\label{eq:ev}  \sum_{i = 1}^{n - 1} (b^{i})^{2} < 1.
 \end{equation}
In normal coordinates,
$$\sum_{i = 1}^{n - 1} (b^{i})^{2}  = |V|^2_{\sigma^*g} = |V|^2_h =
   h(\mu  - h\left(\mu,\partial_\theta)^{-1}\partial_\theta,\mu -  h(\mu,\partial_\theta)^{-1}\partial_\theta\right) 
=1 -2 +  \frac{h(\partial_\theta,\partial_\theta)}{h(\mu,\partial_\theta)^{2}}.$$
Thus (\ref{eq:ev}) holds if
$$ \frac{h(\partial_\theta,\partial_\theta)}{h(\mu,\partial_\theta)^{2}}< 2.$$.
\end{proof}

Using this lemma, we will construct an elliptic PDE with an appropriate inhomogeneous term $ F $ whose solution has 
small $ \calC^{1, \alpha} $-norm. A solution of this PDE plays a crucial role in the proof of Thm.~ \ref{MAIN:thm1}. We begin by finding a function $F$
that is sufficiently large on 
$X\times \{0\}$, 
but with small $ \calL^{p} $-norm on $W$.

\begin{lemma}\label{CG:lemma1} 
Fix  $ p \in \mathbb{N} $ with $ p > n + 1 $ and  a positive $ \delta \ll 1 $. Set $$ C_{g} : = \sup_{X \times \{P\}}  2\left\lvert -{\rm Ric}_{\imath^{*}g}(\mu, \mu)  +  h_{\imath^{*}g}^{2} - \lvert A_{\imath^{*}g} \rvert^{2} \right\rvert. $$
For small enough $\epsilon : = \epsilon(p, C_g, \delta) <1$, there exists 
a positive smooth function $F_{g} = F_{g,\epsilon}:W\to\R$
such that $ F_{g} |_{X \times [-\epsilon/2,\epsilon/2]} 
 = C_{g}$, $ F_{g} \equiv 0 $ outside $ [-\epsilon, \epsilon]$ 
and $ \lVert  F_{g}  \rVert_{\calL^{p}(W, \sigma^{*}g)} < \delta. $
\end{lemma}

Here ${\rm Ric}_{\imath^{*} g}, h_{\imath^{*} g}, A_{\imath^{*} g}$ are the Ricci curvature, mean curvature and second fundamental form of $i^*g.$

\begin{proof} 
Given $ p, C_g, \delta $,  set $ f = C_{g} $ on $ X  = X\times \{0\}$. For a positive constant $ \epsilon \ll 1 $, construct a positive smooth  function $ \phi : [-1, 1] \rightarrow [0,1] $, with $ \phi(t) =  1 $    for $|t|\in \left[0, \frac{\epsilon}{2}\right], $ and  
$ \phi(t) = 0$ for $\lvert t \rvert \in  
\epsilon,1] $.  
Set $F_g = f \cdot \phi :  W \rightarrow \R.$ Then $ F_{g} = C_{g} $ 
on $X \times [-\epsilon/2,\epsilon/2]$ and 
$ F_g>0$ on $W$.  Finally,
$$\Vert F_g\Vert^p_{\calL^{p}(W, \sigma^{*}g)} = \int_{X \times [-\epsilon, \epsilon]} |F_g|^p 
\dvols < \delta,$$
for  $\epsilon$ sufficiently small. 
\end{proof}

\begin{remark}\label{rem:1} Note that $C_g\sim \epsilon^{-1/p}.$  This will be used in the proof of Lem.~\ref{CG:second}. 
\end{remark}

From now on, let $R_\ell$ be the scalar curvature of a metric $\ell$.

Let $\nu = \pm \partial_t$ be the unit outward normal vector field  along $ \partial W $.   
Define the $\calC^{1,\alpha}$ norm of $u:W\to \R$ by choosing a cover $\{(U_i, (x^k_i))\}$ of $W$ and setting
$$\Vert u\Vert_{\calC^{1,\alpha}} = \Vert u\Vert_{\calC^{0,\alpha}} + \sup_{i,k} \Vert \partial_{x_i^k} u\Vert_{\calC^{0,\alpha}}.$$

\begin{proposition}\label{CG:lemma2}
Let $ (W, \sigma^{*} g) $ be as above. 
 Assume that 
 (\ref{eq:intro}) holds. For fixed $ p, \delta, C_{g}, \eta $, there exists a positive 
 $ \epsilon \ll 1 $ and an associated 
$ F_{g} $,  in the sense of Lem.~\ref{CG:lemma1}, such that for  
 fixed $P \in \mathbb{S}^1$,      
\begin{equation}\label{CG:eqn1}
4\nabla_{V} \nabla_{V} u - 4 \Delta_{\sigma^{*} g}  u + R_{g} |_{\sigma(W)}u =
F_{g}  \; {\rm in} \; W, \ u = 0 \; {\rm on} \; \partial W, 
\end{equation}
admits a smooth solution $u:W\to\R$ with
\begin{equation}\label{CG:eqn2}
\lVert u \rVert_{\calC^{1, \alpha}(W)} 
\leqslant \eta.
\end{equation}
\end{proposition}
In the proof and the rest of the paper, norms of gradients like 
$|\nabla_{\sigma^*g} u|_{\sigma^*g}^2$ are usually just denoted by $|\nabla_{\sigma^*g} u|^2.$

\begin{proof}
Denote the $\calL^p$ Sobolev spaces by $\calW^{k,p}.$ Fix 
$p \gg 0$ so that the Sobolev embedding\\
$ \calW^{2, p}(W, \sigma^{*}g) \hookrightarrow \calC^{1, \alpha}(W) $ is a compact inclusion. Define the elliptic operator
\begin{equation*}
L  : = L_{\sigma^*g} = 
4\nabla_{V} \nabla_{V} u - 4 \Delta_{\sigma^{*} g}   + R_{g} |_{\sigma(W)} =4L'  + R_{g} |_{\sigma(W)}
\end{equation*}
on $W$. The Dirichlet boundary condition allows us to rewrite (\ref{CG:eqn1}) on $W$ 
 in the weak form
\begin{align}\label{eq:weak}
\MoveEqLeft{ -\int_{W}  4\nabla_{V} u \cdot \nabla_{V} \phi \  \dvols  + \int_{W} 4 \nabla_{\sigma^{*}g} u \cdot \nabla_{\sigma^{*}g} \phi\  \dvols + \int_{W} u L_{0}(\phi)\  \dvols} \notag\\
& \qquad + \int_{W} R_{g} |_{\sigma(W)} u \phi\ \dvols \\
 &= \int_{W} F_{g} \phi\  \dvols, \forall \phi \in \calC_{c}^{\infty}(W),\notag
\end{align}
where $ L_{0} $ is a first  
order operator determined by $ V $.
Note that 
\begin{equation}\label{R>0} R_h>0 \ \text {on}\  X\times \soo \Rightarrow  R_{g} > 0 \ \text{on}\  M .
\end{equation}
The operator $ L $ is injective: if $ Lu=0$, then combining 
Lem.~\ref{CG:lemma6}, (\ref{R>0}) and the maximum principle gives $ u=0.$ 
It follows from the Fredholm alternative that (\ref{eq:weak}) has a weak solution in $\calW^{1,2}(W) $.  
A standard bootstrapping argument \cite[Thm.~ 1]{Che} then implies $u$ is smooth. By classical elliptic regularity theory on $ (W, \partial W, 
 \sigma^{*}g) $ with 
 Dirichlet boundary condition 
 \cite{Niren4},
 there exists a constant 
$ C_{1}  = C_1(W, \sigma^{*}g,p) $ such that 
\begin{equation}\label{eq:ellest}
    \lVert u \rVert_{\calW^{2, p}(W, \sigma^{*}g)} \leqslant C_{1} \left(\lVert Lu \rVert_{\calL^{p}(W, \sigma^{*}g)} + \lVert u \rVert_{\calL^{p}(W, \sigma^{*}g)} \right),
\end{equation}
for $ u \in \calL^{p}(W, \sigma^{*}g) \cap \calW^{1, 2}(W, \sigma^{*}g) $. 

 We claim that for all $ u \in \calW^{1, 2}(W) \cap \calL^{p}(W) $ solving (\ref{eq:weak}) weakly,
there exists
$ C_{1}'>0 $ such that
\begin{equation}\label{CG:eqn2a}
    \lVert u \rVert_{\calL^{p}(W, \sigma^{*}g)} \leqslant C_{1}' \lVert Lu \rVert_{\calL^{p}(W, \sigma^{*}g)}.
\end{equation}
If not,  there exists a  sequence $ \lbrace u_{n} \rbrace $ solving (\ref{eq:weak}) weakly such that
\begin{equation*}
    \lVert u_{n} \rVert_{\calL^{p}( W, \sigma^{*}g)} = 1, 
     u_n|_{\partial W} = 0 \; {\rm weakly}, \lVert Lu_{n} \rVert_{\calL^{p}(W, \sigma^{*}g)} \leqslant \frac{1}{n}.
\end{equation*}
It follows that
\begin{align*}
    \lVert u_{n} \rVert_{\calW^{2, p}(W, \sigma^{*}g)} & \leqslant C_{1} \left( \lVert  Lu_{n} \rVert_{\calL^{p}( W, \sigma^{*}g)} + \lVert  u_{n} \rVert_{\calL^{p}( W, \sigma^{*}g)} \right) \leqslant 2C_{1}.
\end{align*}
Hence there exists a function $ u \in \calC^{1, \alpha}(W) $ such that
after taking a subsequence, we may assume $ u_{n} \rightarrow u $ in 
$ \calC^{1, \alpha} $. 
Therefore,
\begin{align*}
    \lVert u \rVert_{\calL^{p}(W, \sigma^{*}g)} = 1, & u|_{\partial W} = 0 \;
   {\rm weakly},  
   \\
-\int_{W}  4\nabla_{V} u \cdot \nabla_{V} \phi\  \dvols &+ \int_{W} 4 \nabla_{\sigma^{*}g} u \cdot \nabla_{\sigma^{*}g} \phi\ \dvols + \int_{W} u L_{0}(\phi) \dvols
  \\   
&  + \int_{W} R_{g} |_{\sigma(W)} u \phi\ \dvols
  = 0, \forall \phi \in \calC_{c}^{\infty}(W). 
\end{align*}
Thus $ Lu = 0 $ weakly, so by the injectivity of $L$, we get  $ u = 0 $. This contradicts 
$ \lVert u_n \rVert_{\calL^{p}(W, \sigma^{*}g)} = 1 $. 

By (\ref{eq:ellest}), (\ref{CG:eqn2a}),
\begin{equation}\label{CG:rev2}
\lVert u \rVert_{\calW^{2, p}(W, \sigma^{*}g)}  \leqslant C_{1}(1 + C_{1}') \lVert  Lu \rVert_{\calL^{p}(W, \sigma^{*}g)}.
\end{equation}
Since $u$ is smooth, the Sobolev embedding gives
\begin{equation*}\label{CG:eqn3}
\lVert u \rVert_{\calC^{1, \alpha}(W)} \leqslant C_{2} \lVert u \rVert_{\calW^{2, p}
(W, \sigma^{*}g)} 
\leqslant C_{1}C_{2}(1 + C_{1}') \lVert F_{g} \rVert_{\calL^{p}(W, \sigma^{*}g)},
\end{equation*}
with $ C_{2} = C_2(W, \sigma^*g,L).$ 
Given $ \eta > 0 $, choose $\delta$ such that 
$C_{1}C_{2}(1 + C_{1}') \delta \leqslant \eta$. 
By choosing $\epsilon \ll 1$, we can construct $F_g$
in Lem.~\ref{CG:lemma1} 
 such that $\lVert F_g \rVert_{\calL^{p}(W, \sigma^{*}g)}<\delta.$  
\end{proof}

We need a technical lemma to control the second derivative of the solution $u$ of (\ref{CG:eqn1})
in the $t$ direction (see (\ref{eq:22})). 
\begin{lemma}\label{CG:second}  For any $\eta'>0$,  
we can shrink $\eta$ and $ \epsilon $  in the notation of 
Prop.~\ref{CG:lemma2} 
so that the solution $u$ of (\ref{CG:eqn1}) satisfies
\begin{equation}\label{CG:2nd1}
    \left\lvert \frac{\partial^{2}u}{\partial t^{2}} \right\rvert < \eta'
\end{equation}
on $X \times \left(-\frac{\epsilon}{4},\frac{\epsilon}{4}\right).$  In particular, under the identification $W \approx W\times \{P\}$, we have this estimate on $X\times \{0\} \times \{P\}\subset M.$
\end{lemma}

\begin{proof} 
By the definitions of $ \sigma $ and  
$ g $, the induced metric $ \sigma^{*} g$ and the vector field $ V $ do not contain  $t$ derivatives. In particular,  
$ R_{g} $ is constant on all $ t $-fibers of $M$.  
For fixed positive $\epsilon \ll 1$, set 
\begin{align*}
    U_{1} & := X \times \left(-\frac{\epsilon}{4}, \frac{\epsilon}{4} \right) \subset U_{2} := X \times \left(-\frac{\epsilon}{2}, \frac{\epsilon}{2} \right) \subset U_{3} := X \times \left(-\epsilon, \epsilon \right) 
\subset W, \\
V_{1} & := X \times \left(-\frac{1}{4}, \frac{1}{4} \right) \subset V_{2} := X \times \left(-\frac{1}{2}, \frac{1}{2} \right).
\end{align*}
Applying $ \frac{\partial^{2}}{\partial t^{2}} $ to both sides of (\ref{CG:eqn1}),
we have
\begin{equation}\label{CG:2nd2}
    4\nabla_{V} \nabla_{V} \frac{\partial^{2} u}{\partial t^{2}} - 4\Delta_{\sigma^{*} g} \frac{\partial^{2} u}{\partial t^{2}} + R_{g} |_{\sigma(W)} \frac{\partial^{2} u}{\partial t^{2}} = 0 \Rightarrow  
   L\left(\frac{\partial^{2} u}{\partial t^{2}} \right) = 0 \; {\rm on} \; U_{1},
\end{equation}
since $F_{g} $ is independent of $t$ on $U_1$ by Lem.~\ref{CG:lemma1}. 

Fix $ s > 0 $ such that $ s - \frac{n}{2}  
\geqslant 2 - \frac{n}{p} $, where $ p $ is given in Prop.~\ref{CG:lemma2}. For the metric $ \sigma^{*} g $, we have the $ W^{s, \frac{2(n + 1)}{n - 1}} $-type elliptic estimate for $L $:
\begin{equation}\label{eq:bee}
\lVert v \rVert_{W^{s, \frac{2(n + 1)}{n - 1}}(U_1, \sigma^{*} g)} \leqslant D_{s} \left( \lVert L
v \rVert_{W^{s - 2, \frac{2(n + 1)}{n - 1}}(U_2, \sigma^{*} g)} + \lVert v \rVert_{\calL^{\frac{2(n + 1)}{n - 1}}(U_2, \sigma^{*} g)} \right), \forall v \in H^{s}(U_2, \sigma^{*} g),
\end{equation}
where
$ D_{s} = D_s(L, \sigma^*g, U_1, U_2) $ implicitly depends on $\epsilon.$

The problem is that $D_s$ may blow up
as $\epsilon\to 0$.  
As a result, we need to scale both the time direction $t$ and 
the metric on $ X $
 so that (i) $U_1, U_2$ become $V_1, V_2$, resp., and (ii) the set of solutions of (\ref{CG:2nd2}) is unchanged (up to a time rescaling).
In particular, we set $t' = \epsilon^{-1}t$, $g' = \epsilon^{-2}g $. 

 Since $ g = h + dt^{2} $, the scaling of $g$ forces the scaling
$ h \mapsto h' = \epsilon^{-2}h $. Setting $\psi:W\to W, \psi(x,t') = (x,\epsilon t') = (x,t)$, we have
$$ \psi^*\sigma^*g = \psi^{*}\sigma^{*}(h + dt^{2}) \mapsto \psi^{*}\sigma^{*}(h' + dt^{2}) = \sigma^*(h' + (dt')^2).$$
(\ref{CG:eqn1}) becomes
\begin{equation*}
    4\nabla_{V'} \nabla_{V'} u(\cdot, \epsilon t') - 4\Delta_{\sigma^{*}g'} u(\cdot, \epsilon t') + R_{g'} |_{\sigma(X \times [-\epsilon^{-1}, \epsilon^{-1}])} u(\cdot, \epsilon t') = F_{g'}(\epsilon t') 
\end{equation*}
on  $$ W_\epsilon := \overline{U_3} = X \times [-\epsilon^{-1}, \epsilon^{-1}]_{t'},$$
where we recall the $R_g = R_h$, and $V' = \epsilon V$ is defined by (\ref{eq:Vdef}) for $g'.$  
Taking $ \frac{\partial^{2}}{\partial (t')^{2}} $ on both sides of 
this equation transforms (\ref{CG:2nd2}) into
\begin{align}\label{CG:2nd3}
 & 4\nabla_{V'} \nabla_{V'} \left(
 \frac{\partial^{2}u}{\partial (t')^{2}} \right) - 4\Delta_{\sigma^{*}g'} \left(\frac{\partial^{2}u}{\partial (t')^{2}} \right) + R_{g'} |_{\sigma(W)} \left(\frac{\partial^{2}u}{\partial (t')^{2}} \right) = 0
\Rightarrow L_{\sigma^{*}g'} \left(\frac{\partial^{2}u}{\partial (t')^{2}} \right) = 0 \; {\rm on} \; V_{1}.
\end{align}
Combining (\ref{eq:bee}) for $\sigma^*g'$ with (\ref{CG:2nd3}), we obtain
\begin{equation}\label{eq:1a}
\left\lVert \frac{\partial^{2}u}{\partial (t')^{2}} \right\rVert_{W^{s, \frac{2(n + 1)}{n - 1}}(V_{1}, \sigma^{*}g')}
\leqslant D_{s} \left\lVert \frac{\partial^{2}u}{\partial (t')^{2}} \right\rVert_{\calL^{\frac{2(n + 1)}{n - 1}}(V_{2}, \sigma^{*}g')},
\end{equation}
where $D_s$ and all ``$D$" constants below are independent of 
$\epsilon$ or $u.$  Here the term associated with $L_{\sigma^*g'}v$ in (\ref{eq:bee}) scales by an explicit power of $ \epsilon $, but this term vanishes by (\ref{CG:2nd3}).

Using the basic elliptic estimate, we have
\begin{align}\label{CG:2nd4a}
\left\lVert \frac{\partial^{2}u}{\partial (t')^{2}} \right\rVert_{\calL^{\frac{2(n + 1)}{n- 1}}(V_{2}, \sigma^{*}g' )} & = \epsilon^{2-\frac{n(n - 1)}{2(n + 1)} - \frac{n - 1}{2(n + 1)}} \left\lVert \frac{\partial^{2}u}{\partial t^{2}} \right\rVert_{\calL^{\frac{2(n + 1)}{n- 1}}(U_{2}, \sigma^{*} g)} \epsilon^{2-\frac{n(n - 1)}{2(n + 1)} - \frac{(n - 1)}{2(n + 1)}} \left\lVert u \right\rVert_{W^{2, \frac{2(n + 1)}{n- 1}}(U_{2}, \sigma^{*} g)} \nonumber \\
& \leqslant \epsilon^{2-\frac{n(n - 1)}{2(n + 1)}- \frac{(n - 1)}{2(n + 1)}} D_{1} \left(\lVert 
L u \rVert_{\calL^{\frac{2(n + 1)}{n- 1}}(U_{3}, \sigma^{*}g)} + \left\lVert u \right\rVert_{\calL^{\frac{2(n + 1)}{n- 1}}(U_{3}, \sigma^{*} g)} \right)  \\
& \leqslant \epsilon^{2-\frac{n(n - 1)}{2(n + 1)}- \frac{(n - 1)}{2(n + 1)}} D_{1} \left( \lVert F_{g} \rVert_{\calL^{\frac{2(n + 1)}{n- 1}}(U_{3}, \sigma^{*}g)} + \lVert u \rVert_{\calL^{\frac{2(n + 1)}{n - 1}}(U_{3}, \sigma^{*}g)} \right) \nonumber \\
& \leqslant D_{1} \epsilon^{2- \frac{n(n -1)}{2(n + 1)}} 
\lVert F_{g}(\cdot, \epsilon t') \rVert_{\calL^{\frac{2(n + 1)}{n- 1}}(W_{\epsilon}, \sigma^{*}g))} \nonumber \\
& \qquad + D_{1}\epsilon^{2} \lVert u(\cdot, \epsilon t') \rVert_{\calL^{\frac{2(n + 1)}{n - 1}}(W_1, \sigma^{*}g')}, \nonumber 
\end{align}
where $W_1 = X \times [-1,1]_{t'}.$

Since $ g = h \oplus dt^{2} $, 
the mean curvature for $(M, \partial M )$ satisfies $ h_{g} \equiv 0 $. Therefore, since $ R_{g} > 0 $ and $ \dim M = n + 1 $, the 
relative Yamabe constant \cite{AB2} $\lambda : =  \lambda(M, \partial M, [g])$ 
 for the conformal class $ [g] $ on $ (M, \partial M) $ is positive and satisfies
\begin{equation}\label{eq:ryi}
\lVert v \rVert_{\calL^{\frac{2(n + 1)}{n - 1}}(M, g)}^{2} \leqslant \lambda^{-1}
\left( \frac{4n}{n - 1} \lVert \nabla_{g} v \rVert_{\calL^{2}(M, g)}^{2} + \int_{M} R_{g} v^{2} {\rm dvol}_g \right), \forall v \in H^{1}(M, g).
\end{equation}  
$\lambda$ is a conformal invariant, so we can replace $g$ by $g'$
in (\ref{eq:ryi}).

Let $\bar u = u\circ \pi:M\to\R$ be the pullback of $u:W\to\R$ under the projection $\pi:M\to W$, so $\bar u$ is constant on the $t$-fibers.  Then for $\theta' = \epsilon^{-1}\theta,$
\begin{align}\label{eq:1e}
 \lVert u \rVert_{\calL^{\frac{2(n + 1)}{n - 1}}(W, \sigma^{*} g')}
        &  = \epsilon^{\frac{n -1}{2(n + 1)}} \left(\frac{1}{2\pi} \int_{\mathbb{S}^{1}} d\theta' \int_{W} \lvert u \rvert^{\frac{2(n + 1)}{n - 1}} \text{dvol}_{\sigma^{*} g'}  \right)^{\frac{n - 1}{2(n + 1)}}\\
        &\leqslant \epsilon^{\frac{n -1}{2(n + 1)}} D_{2}  \left( \int_{M}  \lvert 
     \bar{u} \rvert^{\frac{2(n + 1)}{n - 1}} d\text{Vol}_{g'} \right)^{\frac{n - 1}{2(n + 1)}}, \notag
\end{align}
where $D_2$ satisfies $ \text{dvol}_{d(\theta')^2 \oplus \sigma^*g'} \leq D_2^{\frac{2(n+1)}{n-1}}\ \text{dvol}_{g'}$. The fact that $D_2$ is independent of $\epsilon$ (and $u$) 
is the reason we scale $\theta.$

Using (\ref{eq:ryi}) for $g'$, (\ref{eq:1e}) becomes
\begin{align}\label{eq:1b}
\MoveEqLeft{\lVert u \rVert_{\calL^{\frac{2(n + 1)}{n - 1}}(W, \sigma^{*} g')}}\notag \\
& \leqslant \epsilon^{\frac{n -1}{2(n + 1)}} D_{2} 
\lambda^{-1/2}
\left( \frac{4n}{n - 1} \lVert \nabla_{g'} \bar{u} \rVert^2_{\calL^{2}(M, g')} + \int_{M} R_{g'} \bar{u}^{2} \text{dvol}_{g'} \right)^{\frac{1}{2}}\notag \\
& \leqslant \epsilon^{\frac{n -1}{2(n + 1)}} 
D_{2} D_{3} 
\lambda^{-1/2}\left( \frac{4n}{n - 1}\cdot\frac{1}{2\pi}\int_{\mathbb{S}^{1}} d\theta' \int_{W} \lvert \nabla_{\sigma^{*}g'}
u \rvert^{2}_{\sigma^*g'} \text{dvol}_{\sigma^{*}g'} \right.\\
& \qquad \left. + \frac{1}{2\pi}\int_{\mathbb{S}^{1}} d\theta' \int_{W} R_{g'} 
u^{2} \text{dvol}_{\sigma^{*}g'} \right)^{\frac{1}{2}}\notag \\
& = D_{2} D_{3} \lambda^{-1/2}  \epsilon^{-\frac{1}{n + 1}}
\left( \frac{4n}{n - 1} \lVert \nabla_{\sigma^{*} g'} u\rVert_{\calL^{2}(W, \sigma^{*}g')} 
^{2} + \epsilon^{2} \max_{W} \lvert R_{g} \rvert \lVert u \rVert_{\calL^{2}(W, \sigma^{*}g')} 
\right)^{\frac{1}{2}}, \notag
\end{align}
where we have $ \text{dvol}_{g'} \leq D_3^{2}\text{dvol}_{d(\theta')^2 \oplus \sigma^*g'} $. 
We emphasize that the introduction of $\bar u$ lets us use the positivity of $R_g$ and hence $\lambda$ on $M$ to get the estimate (\ref{eq:1b}) for any function $u$ on $W$.

We now estimate $ \lVert F_{g} \rVert_{\calL^{2}\frac{2(n + 1)}{n - 1}(W_{\epsilon}, \sigma^{*}g')}$  by a similar approach. Here we consider $F_g = F_{g}(\cdot, \epsilon t')$ as a smooth function on $W_\epsilon.$ 

First, we observe that for the product metric $ g = h + dt^{2} $ on 
$ M_{a} : = X \times \mathbb{S}^{1} \times [-a, a],$ 
the relative Yamabe constant $ \lambda_{a} : = \lambda(M_{a}, \partial M_{a}, [g]) $ 
is proportional to the relative Yamabe constant $\lambda = \lambda_1$ on $ M = M_{1} $:
\begin{equation}\label{CG:3rd1}
    \lambda_{a} \geqslant a^{\frac{2}{n + 1}}\lambda.
\end{equation}
To see (\ref{CG:3rd1}), we observe that the Yamabe quotient for $ M_{a} $ satisfies
\begin{align*}
    Q_{a}(v) & = \frac{\frac{4n}{n - 1} \lVert \nabla_{g} v \rVert_{\calL^{2}(M_{a}, g)} + \int_{M_{a}} R_{g} v^{2} \dvol}{\lVert v \rVert_{\calL^{\frac{2(n + 1)}{n - 1}}(M_{a}, g)}^{2}} \\
    & = \frac{\frac{4n}{n -1} \int_{-a}^{a} \int_{X \times \mathbb{S}^{1}} \lvert \nabla_{g}v \rvert^{2}d\text{vol}_{h} dt +\int_{-a}^{a} \int_{X \times \mathbb{S}^{1}} R_{g} v^{2}d\text{vol}_{h} dt}{\left( \int_{-a}^{a} \int_{X \times \mathbb{S}^{1}} v^{\frac{2(n + 1)}{n - 1}} d\text{vol}_{h} dt \right)^{\frac{n - 1}{n + 1}}} \\
    & \overset{t = as}{=}  \frac{a \left( \frac{4n}{n -1} \int_{-1}^{1} \int_{X \times \mathbb{S}^{1}} \lvert \nabla_{g}v \rvert^{2}d\text{vol}_{h} ds +\int_{-1}^{1} \int_{X \times \mathbb{S}^{1}} R_{g} v^{2}d\text{vol}_{h} ds\right)}{a^{\frac{n - 1}{n + 1}} \left( \int_{-1}^{1} \int_{X \times \mathbb{S}^{1}} v^{\frac{2(n + 1)}{n - 1}} d\text{vol}_{h} ds \right)^{\frac{n - 1}{n + 1}}} \\
   & \geqslant a^{\frac{2}{n +1}} \lambda.
\end{align*}
We obtain (\ref{CG:3rd1}) by minimizing $ Q_{a}(v) $ over $ v \in H^{1}(M_{a}, g), v \not\equiv 0 $. 

Setting $ \bar{F}_{g} = \pi^{*} F_{g} $ and computing as in 
(\ref{eq:1e}) and (\ref{eq:1b}), we obtain
\begin{equation*}
    \lVert F_{g} \rVert_{\calL^{\frac{2(n + 1)}{n - 1}}(W_{\epsilon}, \sigma^{*}g)} \leqslant \left( \frac{1}{2\pi} \int_{\mathbb{S}^{1}} d\theta \int_{W_{\epsilon}} \lvert F_g \rvert^{\frac{2(n + 1)}{n - 1}} d\text{vol}_{\sigma^{*}g} \right)^{\frac{n -1}{2(n + 1)}} \leqslant D_{2} \left( \int_{M_{\epsilon}} \lvert \bar{F}_{g} \rvert^{\frac{2(n + 1)}{n - 1}} \dvol \right)^{\frac{n -1}{2(n +1)}}.
\end{equation*}
Note that $ \text{supp}(\bar{F}_{g}(\cdot, \epsilon t')) \subset X \times \mathbb{S}^{1} \times [-1, 1]_{t'}$, and that the $ t' $-derivative of $\bar F_g$
 is $ C_{g} O(1) $.  By  (\ref{eq:ryi}) for $ \lambda_{\epsilon^{-1}} $ on $ M_{\epsilon} $,
\begin{align}\label{CG:3rd2}
    \lVert F_{g} \rVert_{\calL^{\frac{2(n + 1)}{n - 1}}(W_{\epsilon}, \sigma^{*}g)} & \leqslant D_{2} \lambda_{\epsilon^{-1}}^{-\frac{1}{2}} \left( \frac{4n}{n - 1} \lVert \nabla_{g} \bar{F}_{g} \rVert_{\calL^{2}(M_{\epsilon}, g)} + \int_{M_{\epsilon}} R_{g} \bar{F}_{g}^{2} \dvol  \right)^{\frac{1}{2}}  \nonumber \\
    & \leqslant D_{2} \epsilon^{\frac{1}{n + 1}} \lambda^{-\frac{1}{2}} \left( \frac{4n}{n - 1} \lVert \nabla_{g} \bar{F}_{g} \rVert_{\calL^{2}(X \times \mathbb{S}^{1} \times [-1, 1]_{t}, g)} + \int_{X \times \mathbb{S}^{1} \times [-1, 1]_{t}} R_{g} \bar{F}_{g}^{2} \dvol  \right)^{\frac{1}{2}}  \\
    & \leqslant D_{2} D_{g} \lambda^{-\frac{1}{2}} \epsilon^{\frac{1}{n + 1}} \cdot C_{g},
    \nonumber
\end{align}
for some $D_g$ independent of $\epsilon.$
We also have
\begin{equation}\label{eq:1c}\left\lVert \frac{\partial^{2}u}{\partial (t')^{2}} \right\rVert_{\calC^{0}(V_{1})}
     \leqslant D_{0} \left\lVert \frac{\partial^{2}u}{\partial (t')^{2}} \right\rVert_{W^{s, \frac{2(n + 1)}{n - 1}}}(V_{1}, \sigma^{*}g) = D_{0} 
     \epsilon^{\frac{n(n - 1)}{2(n + 1)}} \left\lVert \frac{\partial^{2}u}{\partial (t')^{2}} \right\rVert_{W^{s, \frac{2(n + 1)}{n - 1}}}(V_{1}, \sigma^{*}g').
     \end{equation}
The first inequality uses the Sobolev embedding theorem.
For the second inequality, the zeroth order contribution to the $W^{s, \frac{2(n + 1)}{n - 1}}$ norm 
contributes $\epsilon^{\frac{n(n - 1)}{2(n + 1)}}$ from the scaling of the volume form; higher order derivatives in $t$ contribute higher powers of $\epsilon,$ which can be ignored.

This gives
\begin{align*}
\MoveEqLeft{
   \left\lVert \frac{\partial^{2}u}{\partial (t')^{2}} \right\rVert_{\calC^{0}(V_{1})} 
    \stackrel{(18), (25)}{\leqslant} D_{0} D_{s}\epsilon^{\frac{n(n - 1)}{2(n + 1)}} 
    \left\lVert \frac{\partial^{2}u}{\partial (t')^{2}} 
    \right\rVert_{\calL^{\frac{2(n + 1)} 
    {n - 1}}}(V_{2}, \sigma^{*}g'} \\
& \stackrel{(19)}{\leqslant} D_{0} D_{1} D_{s} \epsilon^{2  + \frac{n(n -1)}{2(n + 1)}}\lVert u(\cdot, \epsilon t') \rVert_{\calL^{\frac{2(n + 1)}{n - 1}}(W, \sigma^{*} g')} + D_{0}D_{1}D_{s} \epsilon^{2} \lVert F_{g}(\cdot, \epsilon t') \rVert_{\calL^{\frac{2(n + 1)}{n - 1}}(W_{\epsilon}, \sigma^{*}g')} \\
& \stackrel{(22)}{\leqslant} D_{0} D_{1} D_{2} D_{3} D_{s} \lambda^{-1/2}  \epsilon^{ 2  + \frac{n(n -1)}{2(n + 1)} - \frac{1}{n + 1}} 
\left( \frac{4n}{n - 1} \lVert \nabla_{\sigma^{*} g'}u \rVert_{\calL^{2}(W, \sigma^{*}g')}^{2} + \epsilon^{2} \max_{W} \lvert R_{g} \rvert \lVert u \rVert_{\calL^{2}(W, \sigma^{*}g')} \right)^{\frac{1}{2}} \\
& \qquad +  D_{0}D_{1}D_{s} \epsilon^{2} \lVert F_{g}(\cdot, \epsilon t') \rVert_{\calL^{\frac{2(n + 1)}{n - 1}}(W_{\epsilon}, \sigma^{*}g')} \\
& = D_{0} D_{1} D_{2} D_{3} D_{s} \lambda^{-1/2}  \epsilon^{2} \left( \frac{4n}{n - 1} \lVert \nabla_{\sigma^{*}g} u\rVert_{\calL^{2}(W, \sigma^{*}g)}^{2} + \max_{W} \lvert R_{g} \rvert \lVert u \rVert_{\calL^{2}(W, \sigma^{*}g)} \right)^{\frac{1}{2}} \\
& \qquad + D_{0}D_{1}D_{s} \epsilon^{2} \lVert F_{g}(\cdot, \epsilon t') \rVert_{\calL^{\frac{2(n + 1)}{n - 1}}(W_{\epsilon}, \sigma^{*}g')} \\
&  
\stackrel{(24)}{\leqslant} \bar{D} \epsilon^{2} \eta + \bar{D}' \epsilon^{2} C_{g} \epsilon^{\frac{1}{n + 1}},  
\end{align*}     
\noindent for positive constants $ \bar{D}, \bar{D}' $. 
In the last line, $\eta$ in Prop.~\ref{CG:lemma2} has $C_g$
fixed and $ \epsilon $  to be determined; see Lem.~\ref{CG:lemma1}. By Rem.~\ref{rem:1}, 
$ C_{g} \sim \epsilon^{-\frac{1}{p}} $. Since $ p > n + 1 $, it follows that $ C_{g} \epsilon^{\frac{1}{n + 1}} \ll 1 $ when $ \epsilon $ is sufficiently small. 
If the estimate in  
Lem.~~\ref{CG:lemma1} holds for a fixed $C_{g} $ and some $\epsilon_0$, then it holds for all $\epsilon < \epsilon_0.$ Therefore,
\begin{equation*}
\left\lVert \frac{\partial^{2}u}{\partial t^{2}} \right\rVert_{\calC^{0}(U_{1})} = \epsilon^{-2} \left\lVert \frac{\partial^{2}u}{\partial (t')^{2}} \right\rVert_{\calC^{0}(V_{1})} < \bar{D} \eta + \bar{D}' 
C_g \epsilon^{\frac{1}{n + 1}} < \eta',
\end{equation*}
by  shrinking $\epsilon$ and hence $\eta$ in Prop.~\ref{CG:lemma2} as needed.
\end{proof}

\begin{remark}\label{CG:re1}
(i) The only places in the paper we use the existence of a PSC metric $h$ on $X\times \soo$  are in (\ref{R>0}) and in the positivity of the relative Yamabe invariant above (\ref{eq:ryi}). The only place we use that $X$ and hence $ W $ are oriented is doing the integration by parts/Stokes' Theorem in (\ref{eq:weak}).

(ii) The solution $ u $ of (\ref{CG:eqn1}) on $ W $ may not be positive. However,  $ u_W : = u + 1  $ is positive, provided $ \eta \ll 1 $. 
Let $u_M$ be  the pullback of $u_W$ to $M$ for the projection $M\to W$, and 
let $ u_Y$ be the restriction of $ u_M $ to  
$Y =  X \times \{0\}\times \mathbb{S}^{1} = X\times S^1.$  
Note that 
$u_W$ satisfies
\begin{align}\label{uw}
    4\nabla_{V}\nabla_{V} u_W - 4 \Delta_{\imath^{*} g} u_W + R_{g} |_{\sigma(W)} u_W & = F_g +  R_{g} |_{\sigma(W)} \; {\rm in} \; W, u_W = 1 \; {\rm on} \; \partial W. \\
    \left\lvert \frac{\partial^{2} u_{W}}{\partial t^{2}} \right\rvert < \eta' \; {\rm on} \; X \times \left(-\frac{\epsilon}{4}, \frac{\epsilon}{4}\right) \subset W. \notag
\end{align}

(iii) We use $ e^{2\phi_M} $ as the conformal factor on $ M $, where 
\begin{equation}\label{eq:1f}
    e^{2 \phi_M} := u_M^{\frac{4}{n - 2}}, \tilde g :=  e^{2 \phi_M}g = u_M^{\frac{4}{n - 2}}g \; {\rm on} \; M.
\end{equation} 
We set   $ e^{2 \phi_Y} :
 = u_Y^{\frac{4}{n - 2}} $ as the conformal factors on $ Y = 
  X\times \soo$. (The exponent of $u_M$ is not the usual 
 one. For $X\times \soo$, the exponent is standard since $ n = \dim (X \times \mathbb{S}^{1}) 
 $.) Similarly, set $ e^{2 \phi_W} : 
 = (u_W)^{\frac{4}{n - 2}} $ on $W$.
As pullbacks,  $ \phi_M $ and $ u_M $ are constant
 in the $ \mathbb{S}^{1} $-direction:
\begin{equation*}
    \frac{\partial \phi_M}{\partial \theta} = \frac{\partial u_M}{\partial \theta} = 0.
\end{equation*}

(iv) The PDE (\ref{CG:eqn1}) is related to the conformal change 
\begin{equation}\label{eq:uphi} \tilde{g} = u^{\frac{4}{n - 2}} g,
\end{equation}
for a positive, smooth function $ u $ and a metric $g$ on any manifold. For $\phi$ defined by 
\begin{equation}\label{eq:uphi2}e^{2\phi} = u^{\frac{4}{n - 2}},
\end{equation}
it is straightforward to check by differentiating (\ref{eq:uphi2}) that the Laplacians of these two conformal factors are related by
\begin{equation*}
    u^{-\frac{n+2}{n- 2}}\left(-\frac{4}{n - 2} \Delta_{g} u 
    \right) = e^{-2\phi} \left(  
    - 2 \Delta_{g} \phi - (n - 2) \lvert \nabla_{g} \phi \rvert^{2} \right), n \geqslant 3,
\end{equation*}
and their gradients are related by
\begin{equation*}
  e^{-2\phi} \nabla_{g} \phi = \frac{2}{n - 2} u^{-\frac{n+2}{n-2}} \nabla_{g} u,\  e^{-2\phi} \lvert \nabla_{g} \phi \rvert^{2} = \left( \frac{2}{n - 2} \right)^{2}u^{-\frac{2n}{n - 2}} \lvert \nabla_{g} u \rvert^{2}. 
\end{equation*} 
Finally, we will use the following relation between $u_Y$ and $\phi_Y$, again given by differentiating
(\ref{eq:uphi2}):
$$e^{-2\phi_{Y}} \left( 2(n - 2) \nabla_{V} \nabla_{V} \phi_Y + (n - 2)^{2} \nabla_{V} \phi_Y \nabla_{V} \phi_Y \right)   
= 4u_{Y}^{-\frac{n + 2}{n - 2}}   \nabla_{V} \nabla_{V} u_{Y}.$$

(v) We have
\begin{align*}u_Y(x,P) = u_M(x,0,P) = u_W(x,0),\\ 
\phi_Y(x,P) = \phi_M(x,0,P) = \phi_W(x,0).
\end{align*}

\end{remark}

\medskip

As mentioned above,
we need to compare the Laplacians on $(M,g)$,  $(W, \sigma^*g)$, and  $ (X \times \mathbb{S}^{1},h) $.
\begin{lemma}\label{CG:lemma3}
In the notation of (\ref{diagram}) and Remark \ref{CG:re1}(iii), we have 
\begin{equation}\label{CG:eqn4}
 \Delta_{g} u_M(w,P)
 - \Delta_{\sigma^{*}g} 
 u_W(w) = B_{1}(u_W, g)(w), 
\end{equation}
for $w\in W$ and fixed $P\in \mathbb{S}^1.$
Here $ B_{1}(u_W, g) = B_{1}(u_W, g,P)$ 
is a globally defined first order operator on $u_W\in C^\infty(W)$; in particular, the left hand side of (\ref{CG:eqn4}) has no derivatives in the $ \mathbb{S}^{1} $-direction.
Similarly, we have
\begin{equation}\label{CG:eqn4a}
\Delta_{g} u_M (y, 0) - \Delta_{\imath^{*}g} u_Y(y) = 
\frac{\partial^{2} u_M}{\partial t^{2}}(y,0),
\end{equation}
for $y \in X\times \mathbb{S}^{1} $.
\end{lemma}

\begin{proof}
In (\ref{CG:eqn4}), for $\pi:M\to W$, we are computing
\begin{equation}\label{CG:LB1}
\Delta_{g} u_{M} - \pi^{*}(\Delta_{\sigma^{*}g}u_{W}) = \Delta_{g} (u_{W} \circ \pi) - \pi^{*}(\Delta_{\sigma^{*}g}(u_{W})) = \Delta_{g} (u_{W} \circ \pi) - \pi^{*}(\Delta_{\sigma^{*}g}(u_{M} \circ \sigma))
\end{equation}
in a tubular neighborhood $T$ of the hypersurface $ H := W \times \lbrace P \rbrace \subset M$, 
and then restricting
to $H$. 
By taking the level surfaces for the $g$-distance function to $H$, we can take a chart $(x^1,\ldots,x^n, t)$ in $M$ containing a chart $(x^1,\ldots,x^{n - 1}, t)$ for  
$W\times \{P\}$, with $ t : = x^{n + 1} $  the coordinate on $ [-1, 1] 
$ and $ \partial_{ x^{n}}$ normal to the level surfaces.
If $x^n(P)=0$ in  these coordinates, we can write $ h = h_{X, x_{n}} \oplus d(x^{n})^2 $ 
and $ g = h_{X, x_{n}} \oplus d(x^{n})^2 \oplus dt^{2} $ near  
$ W \times \lbrace P \rbrace $. 
Therefore, $ \sigma^{*}g = h_{X, 0} \oplus dt^{2} $ on $ W \approx W \times \lbrace P \rbrace $, 
and $ g^{ij} = (\sigma^{*}g)^{ij}, i, j = 1, \dotso, n - 1 $ on $ W $.
In summary, on $T$,
\begin{equation*}
    (g_{ij}) = \begin{pmatrix} g_{11} & \dotso & g_{1, n - 1} & 0 & 0 \\
    \vdots & \dotso & \vdots & \vdots & \vdots \\
    g_{n -1, 1} & \dotso & g_{n -1, n - 1} & 0 & 0 \\
    0 & \dotso & 0 & 1 & 0 \\
    0 & \dotso & 0 & 0 & 1 \end{pmatrix},
\end{equation*}
with $ (g^{ij}) $ the inverse of the whole matrix, and $ (\sigma^{*}g)^{ij} $  the inverse of the minor given by deleting the $n^{\rm th}$ 
row and 
column.
Then
\begin{equation*}
\Delta_{g} u_M = \sum_{i, j= 1}^{n- 1} g^{ij} \frac{\partial^{2} u_M}{\partial x^{i} \partial x^{j}} + \frac{\partial^{2} u_{M}}{\partial (x^{n})^{2}} + \frac{\partial^{2}u_{M}}{\partial t^{2}} - \sum_{i, j, k = 1}^{n + 1} g^{ij} \Gamma_{ij, g}^{k} \frac{\partial u_M}{\partial x^{k}}.
\end{equation*}

The projection map $ \pi : M \rightarrow W $  
is $ \pi (x^{1}, \dotso, x^{n -1}, \theta, t) = (x^{1}, \dotso, x^{n-1}, t)$, 
so $\frac{\partial u_{M}}{\partial x^{n}} = \frac{\partial (u_{W} \circ \pi)}{\partial x^{n}} =0.$
Similarly, it is easy to check that 
 $$   \frac{\partial u_{W}}{\partial x^{i}}\biggl|_{w} = \frac{\partial (u_{M} \circ \sigma)}{\partial x^{i}}\biggl|_{w} = 
  \frac{\partial u_{M}}{\partial x^{i}}\biggl|_{(w,P)}, i = 1, \dotso, n - 1, n + 1.$$
 On $ W$, 
\begin{equation*}
 \Delta_{\sigma^{*}g} u_W   = \sum_{i, j= 1}^{n-1} (\sigma^{*}g)^{ij} \frac{\partial^{2} 
 u_W}{\partial x^{i} \partial x^{j}} + \frac{\partial^{2}
 u_W}{\partial t^{2}} - (\sigma^{*}g)^{ij} \Gamma_{ij, \sigma^{*}g}^{k} \frac{\partial 
 u_W}{\partial x^{k}}.
\end{equation*}
Since  
$ (g^{ij}) = (\sigma^{*}g)^{ij} $ on $ W $, 
we have
\begin{align*}
   & 
   \Delta_{g}u_{M}(w, P) - \Delta_{\sigma^{*}g}u_{W}(w) \\
   & \qquad = \left( \sum_{i, j= 1}^{n- 1} g^{ij} \frac{\partial^{2} u_M
   }{\partial x^{i} \partial x^{j}}\biggl|_{(w, P)} + \frac{\partial^{2}u_{M}
   }{\partial t^{2}}\biggl|_{(w, P)} - \sum_{i, j,k = 1}^{n + 1} g^{ij} \Gamma_{ij, g}^{k} \frac{\partial u_W
   }{\partial x^{k}}\biggl|_w \right) \\
   & \qquad\qquad - \left( \sum_{i, j
   = 1}^{n-1} (\sigma^{*}g)^{ij} \frac{\partial^{2}u_{M}
   }{\partial x^{i} \partial x^{j}}\biggl|_{(w, P)} + \frac{\partial^{2}u_{M}
   }{\partial t^{2}}\biggl|_{(w, P)} - \sum_{i, j, k \neq n}(\sigma^{*}g)^{ij} \Gamma_{ij, \sigma^{*}g}^{k} \frac{\partial u_{W}
   }{\partial x^{k}}\biggl|_w \right) \\
   & \qquad : = B_{1}(u_{W}, g)(w).
\end{align*}

 (\ref{CG:eqn4a}) follows from the fact that $ g = \imath^{*} g \oplus dt^{2} $ with $ t = x^{n + 1} $, and hence all Christoffel symbols involving the $ x^{n + 1} $ direction vanish.
\end{proof}

\begin{lemma}\label{lem:4.6}
   For fixed $\eta$ 
    in Prop.~\ref{CG:lemma2}, for  
 fixed $P\in \soo$ and $u_W$ as in Remark \ref{CG:re1}(ii), there exists a constant $K_1 = K_1(\eta)>0 $
 such that
    \begin{equation}\label{eq:11b} 
  4 |B_1(u_W,g,P)(w)|
  <K_{1},     
    \end{equation}  
    for all $w\in  
    W.$ We have $K_1\to 0$ as $\eta\to 0.$
\end{lemma}

\begin{proof} 
  $B_1$ is a 
   first order operator in $u$ with smooth coefficients, and which does not contain derivatives of $u$ in the normal or $\soo$ directions. Since the coefficients of $B_1$ are bounded on $W$ and we have $\Vert u\Vert_{\calC^{1,\alpha}}<\eta$, the result follows.
  \end{proof}

\section{Proof of The Main Theorem}
In this section, we prove   
Thm.~ \ref{MAIN:thm1}.
One direction of the $\soo$-Stability Conjecture direction is trivial: if $ X $ admits a PSC metric $ g_{0} $, then the product metric on $ X \times \mathbb{S}^{1} $ has PSC.

As a possible approach to the $\soo$-stability conjecture in the other direction,
if we start with $R_h = R_{\imath^{*}g} >0$ on $ X \times \soo $, then $ R_{\imath^{*}g} |_{X \times \lbrace P \rbrace} >0$. In the notation of (\ref{diagram}), if we somehow knew that 
$  R_{\imath^{*}g} |_{X \times \lbrace P \rbrace} - R_{\tau^{*}\imath^{*}g} |_{X} 
<0$
 by some Gauss-Codazzi manipulations,
then $ R_{\tau^{*}\imath^*g} > 0 $ on $ X$. Unfortunately, it is very hard to estimate the second order derivatives of the metric in the Gauss-Codazzi equation. Instead, we make two modifications: (i) we introduce the extra dimension $[-1,1]$ in order to produce the elliptic operator $L'$ on $W$ in Lem.~\ref{CG:lemma6}; (ii) motivated by \cite{RRX}, we consider a conformal transformation of  $g$ to $\tilde g$ as in (\ref{eq:1f}).  In terms of the diagram (\ref{diagram}), the extra dimension in (i) allows us to transfer geometric information up the right hand side, {\em i.e.,} from $(W, \sigma^*g)$ to 
$(M, g)$. The conformal transformation in (ii) produces 
$(M,\tilde g)$ in the upper right corner, and Gauss-Codazzi moves us first to the upper left corner and then down to
 $(X, \tau^* \imath^* \tilde g)$ in the lower left corner.  In the end, we will prove that 
 $  \tau^* \imath^* \tilde g$ is a PSC metric on $X$, provided 
 (\ref{eq:intro}) holds.

\begin{theorem}\label{MAIN:thm1}
Let $ X $ be an oriented closed manifold with $ \dim X = n - 1 \geqslant 2 $. If $ X \times \mathbb{S}^{1} $ admits a PSC metric $ h $ such that for some $P\in \soo$, we have
\begin{equation}\label{eq:intro1}
\angle_h(\mu, \partial_\theta) < \frac{\pi}{4} 
\end{equation}
on $X\times \{P\},$
then $ X \approx X \times \{P\}$ admits a PSC metric.
\end{theorem}

Most of the proof works with a general conformal transformation $u_W $ of $\sigma^*g$ on $W$; it is only in (\ref{eq:26}) that we specify that $u_W$ satisfies (\ref{CG:eqn1}).

\begin{proof} 
We collect a few facts in preparation for the proof.  
While we choose $\eta$ and $\delta$ sufficiently small in the bullets below, the key point is that all bulleted estimates  hold for $\epsilon$ sufficiently small.
\begin{itemize}
\item{}We fix 
$ p \in \mathbb{N}, p \gg 1 $, 
such that $ \calW^{2, p}(W, \sigma^{*}g) \hookrightarrow \calC^{1, \alpha}(W) $ is a compact inclusion.

\item{}  Pick $\eta \ll 1$ and $\delta = \delta(\eta)>0$ such that for the function $F_g = F(p,C_{g},\delta)$ in Lem.~\ref{CG:lemma1}, the solution $u_W$ of (\ref{uw}) satisfies 
    $\frac{1}{2} < \Vert u_W \Vert_{\calC^{0}} < \frac{3}{2} $; automatically, $\sup_{i,k} \Vert \partial_{x^i_k} u_W\Vert_{C^{0,\alpha}} \leqslant \eta$
    by Prop.~\ref{CG:lemma2} and the notation above the Proposition. 
Set $R_{\rm inf} := \inf_Y R_{h} = \inf_M R_g>0$. We further shrink $\epsilon $
as necessary so that 
\begin{equation}\label{etaprime}\eta'< \frac{1}{16} R_{\rm  inf}
\end{equation} 
in (\ref{CG:2nd1}).

\item{} Fix $P\in \soo.$ Then
 $ |B_1(u_W,g,P)(w)|< K_1= K_1(\eta) $ 
by Lem.~\ref{lem:4.6}. Recall that $K_1\to 0$ as $\eta\to 0$, so 
we can shrink $\eta$ further (by shrinking $\epsilon$) so that
\begin{equation}\label{K1} K_1 < \frac{1}{16} R_{\rm inf}.
\end{equation} 

\item{} By Lem.~\ref{CG:lemma1}, 
$ C_{g} $ satisfies
\begin{align}\label{eq:10a}
C_{g} & = 2 \max_{x\in X } \left\lvert -2  {\rm Ric}_{\imath^{*}g}(\mu, \mu)(x, P)  + h_{\imath^{*} g}^{2}(x,P)  - \lvert A_{\imath^{*} g} \rvert^{2}(x,P) \right\rvert \nonumber \\
& > \frac{3}{2} \left( -2  {\rm Ric}_{\imath^{*}g}(\mu, \mu)(x, P)   + h_{\imath^{*} g}^{2}(x,P)  - \lvert A_{\imath^{*} g} \rvert^{2}(x,P)  \right),\ \forall x \in X. 
\end{align}

\item{} Using the lower bound on $\Vert u_W\Vert$ above and $\Vert 
u\Vert_{\calC^{1,\alpha}} < \eta$ by Prop.~\ref{CG:lemma2}, we can choose $ \eta \ll 1, $ by 
shrinking $ \epsilon $ if necessary, so that for all 
$x\in X,$
$$     K_{2} (u_W, g)(x)  : =  \frac{4}{n - 2}  \left( 
\frac{\lvert \nabla_{\imath^{*}g} u_W(x,0) \rvert^{2}_g}{u_W(x,0)} + n \frac{\nabla_{V} u_W(x,0) \nabla_{V} u_W(x,0) }{u_W(x,0) } \right) $$
has     
 \begin{equation}\label{MAIN:eqn6} \lvert K_{2}(u_{W}, g)(x) \rvert \leqslant C_{3} \eta  < \frac{1}{4}R_{\rm inf},
\end{equation}
where$ C_{3} = C_{3}(g, P, n) $.
This estimate will be used in the last step of the proof.

\item{}  All derivatives
of $u_M, \phi_M, u_Y, \phi_Y$ vanish in the $\soo$ direction, by their definitions  in Remark \ref{CG:re1}(ii),(iii).  In particular, $\nabla_\mu\nabla_\mu \phi_Y = \nabla_V\nabla_V\phi_Y$.

\end{itemize}

We now work on the left side of (\ref{diagram}) using the Gauss-Codazzi equation.  
Gauss-Codazzi and the conformal transformation law for the second fundamental form on 
 $ X \times \lbrace P \rbrace \subset X \times \mathbb{S}^{1} $ give
\begin{align*}
R_{\tau^{*} \imath^{*} \tilde{g}} & = R_{\imath^{*} \tilde{g}} - 2  \text{Ric}_{\imath^{*} \tilde{g}}\left(e^{-\phi_Y} \mu, e^{-\phi_Y} \mu \right) + h_{\imath^{*} \tilde{g}}^{2} - \lvert A_{\imath^{*} \tilde{g}} \rvert^{2}, \\
A_{\imath^{*}\tilde{g}}(V_{1}, V_{2}) & = e^{\phi_Y}A_{\imath^{*}g}(V_{1}, V_{2}) + e^{\phi_Y} \frac{\partial \phi_Y}{\partial \mu} g(V_{1}, V_{2}), \forall V_{1}, V_{2} \in T(X \times \lbrace P \rbrace).
\end{align*}
 It follows that
\begin{align*}
   \lvert A_{\imath^{*} \tilde{g}} \rvert^{2} & = e^{-2\phi_Y} \left( \lvert A_{\imath^{*} g} \rvert^{2} + 
   2 n h_{\imath^{*}g} \frac{\partial \phi_Y}{\partial \mu} +  n^2
   \left(\frac{\partial \phi_Y}{\partial \mu} \right)^{2}  \right),\\
    h_{\imath^{*} \tilde{g}}^{2} & = e^{-2\phi_Y} \left( h_{\imath^{*} g}^{2} +
     2    n\frac{\partial \phi_Y}{\partial \mu} h_{\imath^{*}g} + 
       n^2\left(\frac{\partial \phi_Y}{\partial \mu} \right)^{2} \right).
\end{align*}
The formulas for the conformal transformation of $ \text{Ric}_{\imath^{*} \tilde{g}} $ and 
$R_{\imath^{*} \tilde{g}} $ on $ X \times \mathbb{S}^{1} $ with conformal factor $ e^{2\phi_Y} $ are
\begin{align}\label{eq:ric}
    \text{Ric}_{\imath^{*} \tilde{g}}\left(e^{-\phi_Y} \mu, e^{-\phi_Y} \mu\right) & = e^{-2\phi_Y}\left(\text{Ric}_{\imath^{*}g}\left(\mu, \mu\right) - (n - 2)(\nabla_{\mu} \nabla_{\mu} \phi_Y - \nabla_{\mu} \phi_Y \nabla_{\mu} \phi_Y)\right)\notag \\
    & \qquad - e^{-2\phi_Y}\left(\Delta_{\imath^{*}g} \phi_Y + (n - 2) \lvert \nabla_{\imath^{*} g} \phi_Y \rvert^{2}\right)g_{nn},\\
  R_{\imath^{*} \tilde{g}}  &=  e^{-2\phi_Y} (R_{\imath^*g} -2(n-1)\Delta_{i^*g}\phi_Y
  -(n-1)(n-2)|\nabla_{\imath^*g}\phi_Y|^2).\notag
  \end{align}
Since $ g_{nn} = 1 $, we have  
\begin{align}\label{MAIN:eqn5}
R_{\tau^{*} \imath^{*} \tilde{g}} & = e^{-2\phi_Y} \left(R_{\imath^{*} g} - 2 \text{Ric}_{\imath^{*}g}\left(\mu, \mu\right) + h_{\imath^{*} g}^{2} - \lvert A_{\imath^{*} g} \rvert^{2} \right)\notag \\
    & \qquad + e^{-2\phi_Y} \left( - 2(n - 1) 
    \Delta_{\imath^{*}g} \phi_Y - (n - 1)(n - 2) \lvert \nabla_{\imath^{*}g} \phi_Y \rvert^{2} \right) \\
    &\qquad + e^{-2\phi_Y} \left(2(n - 2) \left(\nabla_{\mu} \nabla_{\mu} \phi_Y - \nabla_{\mu} \phi_Y \nabla_{\mu} \phi_Y \right) + 2 \left(\Delta_{\imath^{*} g} \phi_Y + (n - 2) \lvert \nabla_{\imath^{*} g} \phi_Y \rvert^{2} \right) \right).\notag 
\end{align}
Here $R_{\tau^{*} \imath^{*} \tilde{g}}$ is evaluated at $x\in X$, and the right hand side of (\ref{MAIN:eqn5}) is evaluated at $y=(x,P)\in X\times \{P\}\subset Y.$

Now we replace the second order terms in $\phi_Y$ in  (\ref{MAIN:eqn5}) with terms in $u_Y $. 
Using Remark \ref{CG:re1}(iv), we get
\begin{align}R_{\tau^{*} \imath^{*} \tilde{g}} 
& = u_{Y}^{-\frac{n + 2}{n - 2}} \left( - 2 \text{Ric}_{\imath^{*}g}\left(\mu, \mu\right) u_{Y} + h_{\imath^{*} g}^{2} u_{Y} - \lvert A_{\imath^{*} g} \rvert^{2} u_{Y} \right) \\
      & \qquad + u_{Y}^{-\frac{n + 2}{n - 2}} \left( 4 \nabla_{V} \nabla_{V} u_{Y} - 4 \Delta_{\imath^{*} g} u_{Y} + R_{\imath^{*} g} u_{Y} \right)\notag \\
     & \qquad + e^{-2\phi_{Y}} \left( (n - 2) \lvert \nabla_{\imath^{*}g} \phi_Y \rvert^{2}  - n(n - 2) \nabla_{V} \phi_Y \nabla_{V} \phi_Y \right) \notag.
\end{align}
In the next equation, we use  Remark \ref{CG:re1}(ii), Remark \ref{CG:re1}(iv), Lem.~\ref{CG:second}, and
(\ref{CG:eqn4a}). In the first equality, the left hand side is evaluated at $x\in X$,  and all terms on the right hand side (including $h_{\imath^{*} g}^{2}$ and $\lvert A_{\imath^{*} g} \rvert^{2}$) are evaluated at $y = (x,P)\in X\times \{P\}$; in the next inequality, the term $\Delta_g u_M$  is evaluated at $(x,0,P)\in M$; in the last equality, $K_2$ is evaluated at $x.$
We obtain
\begin{align}\label{eq:22}
     R_{\tau^{*} \imath^{*} \tilde{g}} & 
    =u_Y^{-\frac{n + 2}{n - 2}} \left(- 2 \text{Ric}_{\imath^{*}g}\left(\mu, \mu \right) u_Y + h_{\imath^{*} g}^{2} u_Y - \lvert A_{\imath^{*} g} \rvert^{2} u_Y \right)\notag \\
    & \qquad + u_Y^{-\frac{n + 2}{n - 2}} \left( 4\nabla_{V} u_Y \nabla_{V} u_Y -4 \Delta_{g} u_M +  4\frac{\partial^{2} u_{Y}}{\partial t^{2}} + R_{\imath^{*} g} u_Y \right)\notag \\
    & \qquad  + e^{-2\phi_Y}  \left( -(n - 2)  \lvert \nabla_{\imath^{*}g} \phi_Y \rvert^{2} - n(n - 2)  \nabla_{V} \phi_Y \nabla_{V} \phi_Y \right) \notag\\
   & \geqslant u_Y^{-\frac{n + 2}{n - 2}} \left(- 2 \text{Ric}_{\imath^{*}g}\left(\mu, \mu \right) u_Y + h_{\imath^{*} g}^{2} u_Y - \lvert A_{\imath^{*} g} \rvert^{2} u_Y \right) \\
    & \qquad + u_Y^{-\frac{n + 2}{n - 2}} \left( 4\nabla_{V} u_Y \nabla_{V} u_Y -4 \Delta_{g} u_M 
    - 4\eta' + R_{\imath^{*} g} u_Y \right) \notag\\
    & \qquad  +   u_Y^{-\frac{n + 2}{n - 2}}\cdot \frac{4}{n - 2}  \left( - \frac{\lvert \nabla_{\imath^{*}g} u_Y \rvert^{2}}{u_Y} - n \frac{\nabla_{V} u_Y \nabla_{V} u_Y}{u_Y}\right) \notag\\
    &  = u_Y^{-\frac{n + 2}{n - 2}} \left(- 2 \text{Ric}_{\imath^{*}g}\left(\mu, \mu \right) u_Y + h_{\imath^{*} g}^{2} u_Y - \lvert A_{\imath^{*} g} \rvert^{2} u_Y \right)\notag \\
    & \qquad + u_Y^{-\frac{n + 2}{n - 2}} \left( 4\nabla_{V} u_Y \nabla_{V} u_Y -4 \Delta_{g} u_M 
    + R_{\imath^{*} g} u_Y \right) \notag\\
    & \qquad  + u_Y^{-\frac{n + 2}{n - 2}}\left(-K_{2}(u_{W}, g) - 4\eta' \right). \notag
\end{align}

Now we replace $\Delta_g$ with $\Delta_{\sigma^*g},$ which is a transfer of information from lower right to upper right in (\ref{diagram}).  In the first inequality in (\ref{eq:26}), we use (\ref{CG:eqn4});
in the next equality, we use that $u_W(x,0) = u_Y(x,P)$ on $X\times \{P\}\times \{0\}$ by Remark \ref{CG:re1}(iv), and $ R_{g} = R_{\imath^{*}g} $, since $ g 
$ is a product metric; in the last equality,
 we use (\ref{CG:eqn1}), (\ref{uw}) 
and the fact that $F_g = C_g$ on $X\times \{0\}.$ Again, 
$R_{\tau^{*} \imath^{*} \tilde{g}}$ and and 
$K_2(u_W,g)$ are evaluated at $x$; $u_Y, h_{\imath^*g},$ and  $ |A_{\imath^*g}|^2 $ are evaluated at $y = (x,P)$;  $u_W$ and $ B_1(u_W,g)$ are evaluated at $w = (x,0)\in W;$ $R_g$ is evaluated at $(x,0,P) \in M.$  
Then for $u_W$ solving (\ref{CG:eqn1}), we have
\begin{align}\label{eq:26}
     R_{\tau^{*} \imath^{*} \tilde{g} }& \geqslant u_Y^{-\frac{n + 2}{n - 2}} \left(- 2 \text{Ric}_{\imath^{*}g}\left(\mu, \mu \right) u_Y + h_{\imath^{*} g}^{2} u_Y - \lvert A_{\imath^{*} g} \rvert^{2} u_Y \right) \notag\\
      & \qquad + u_Y^{-\frac{n + 2}{n - 2}} \left( 4 \nabla_{V} \nabla_{V} u_Y -4\Delta_{\sigma^{*}g} u_Y + R_{\imath^{*} g} u_Y \right) \notag\\
      & \qquad - u_Y^{-\frac{n + 2}{n - 2}} \left(4 B_{1}(u_Y, g) \right) + u_Y^{-\frac{n + 2}{n - 2}} \left(-K_{2}(u_{W}, g) - 4 \eta'\right) \notag\\
      & = \left(u_W \right)^{-\frac{n + 2}{n - 2}} \left(- 2 \text{Ric}_{\imath^{*}g}\left(\mu, \mu \right) u_{W} + h_{\imath^{*} g}^{2} u_{W} - \lvert A_{\imath^{*} g} \rvert^{2}
       u_{W} \right)\\
      & \qquad + \left(u_W \right)^{-\frac{n + 2}{n - 2}} \left( 4 \nabla_{V} \nabla_{V} u_{W} -4 \Delta_{\sigma^{*}g} u_{W} + R_{g} |_{\sigma(W)} u_{W} \right)\notag \\
      & \qquad  - \left(u_W \right)^{-\frac{n + 2}{n - 2}} \left(4 B_{1}(u_{W}, g) \right) + \left(u_W \right)^{-\frac{n + 2}{n - 2}} \left(-K_{2}(u_{W}, g) - 4\eta' \right) \notag\\
      & = \left(u_W \right)^{-\frac{n + 2}{n - 2}} \left(- 2 \text{Ric}_{\imath^{*}g}\left(\mu, \mu \right) u_{W} + h_{\imath^{*} g}^{2} u_{W} - \lvert A_{\imath^{*} g} \rvert^{2} u_{W} \right) \notag\\
      & \qquad + \left(u_W \right)^{-\frac{n + 2}{n - 2}} \left( C_{g} + R_{g} -  4 B_{1}(u_{W}, g) - K_{2}(u_{W}, g) - 4\eta' \right).\notag 
\end{align}
It follows from

(\ref{R>0}),  and (\ref{etaprime}) --  
(\ref{MAIN:eqn6}) that
\begin{align*}\label{MAIN:eqn4}
R_{\tau^{*} \imath^{*} \tilde{g}} &> u_W^{-\frac{n + 2}{n - 2}} \left( C_{g} - \frac{3}{2}\max_{x \in X}\left\lvert - 2  \text{Ric}_{\imath^{*}g}\left(\mu, \mu \right)(y, P)  + h_{\imath^{*} g}^{2}(x,P)  - \lvert A_{\imath^{*} g} \rvert^{2}(x,P)  \right\rvert + R_{\rm inf} -  \frac{3}{4} R_{\rm inf}
\right)\\
&> 0.\notag
\end{align*}
\end{proof}

\begin{remark}\label{MAIN:re1}
The technical difficulty in the proof is assuring that the operator $L'$ in Lem.~\ref{CG:lemma6}
is elliptic, so that we get a smooth conformal factor $u_W.$  The difficult term $\nabla_\mu\nabla_\mu$
in $L'$ appears because (i) we have to work on $M = X \times \soo\times [0,1]$, as explained before 
Lem.~\ref{CG:lemma6}, and (ii)  to pass from $M$ to $X$, we have to use Gauss-Codazzi for 
$\tilde g$ twice:
\begin{align*}
    R_{\tau^{*} \imath^{*} \tilde{g}} & = R_{\imath^{*} \tilde{g}} - 2\text{Ric}_{\imath^{*} \tilde{g}}\left(e^{-2\phi_Y} \mu, e^{-2\phi_Y} \mu \right) +  h_{\imath^{*} \tilde{g}}^{2} - \lvert A_{\imath^{*} \tilde{g}} \rvert^{2} \\
    & = R_{\tilde{g}} - 2\text{Ric}_{\tilde{g}}\left(e^{-2 \phi_M} \nu', e^{-2 \phi_M} \nu' \right) + h_{\tilde{g}}^{2} - \lvert A_{\tilde{g}} \rvert^{2} \\
    & \qquad - 2\text{Ric}_{\imath^{*} \tilde{g}}\left(e^{-2\phi_Y} \mu, e^{-2\phi_Y} \mu \right) +  h_{\imath^{*} \tilde{g}}^{2} - \lvert A_{\imath^{*} \tilde{g}} \rvert^{2}.
\end{align*}
The term $\nabla_\mu\nabla_\mu\phi_Y$ then appears in  
(\ref{eq:ric}) from the  
formulas for conformal transformations from  
$\text{Ric}_{\imath^{*} \tilde{g}}(e^{-\phi_Y}\mu,e^{-\phi_Y}\mu)$
 to $\text{Ric}_{\imath^* g}(\mu,\mu)$. We must use these formulas, since the theorem's hypothesis involves the metric $h = \imath^*g.$
\end{remark}

There are several known classes of counterexamples to the $\soo$-stability conjecture,
including  odd degree hypersurfaces in $\C\mathbb{P}^3$ 
\cite[Rmk.~1.25]{JR},
and $M\ \# \ k \overline{\C\mathbb{P}^2}$ for a simply connected K\"ahler surface $M$ and a positive integer $k$
\cite[Rmk.~5]{KS}. These are $4$-manifolds $X$ which admit no PSC metric, but such that $X\times S^1$ admits a PSC metric.   For these and any other counterexamples, the PSC metric on $X\times S^1$ 
has large angles in the following sense:
\begin{corollary}\label{MAIN:cor0} If $X$ admits no PSC metric, but $X\times \soo$ admits a PSC metric $h$, then  for all $P\in \soo$, there exists $x = x(P)$ such that 
\begin{equation*}
\angle_h(\mu, \partial_\theta))_{(x,P)} \geq \frac{\pi}{4}.
\end{equation*}
\end{corollary}

There is an analytic consequence of Thm.~ \ref{MAIN:thm1}. 
Let $ \lambda_{g}(X) $ be the Yamabe invariant of the closed manifold $ X $ with respect to the conformal class $ [g] $. 
\begin{corollary}\label{MAIN:cor1}
Let $ X \times \mathbb{S}^{1} $ 
be an oriented closed manifold with metric $h$
satisfying (\ref{eq:intro}) on some $X\times \{P\}.$  
Then
\begin{equation*}
    \lambda_{h}(X \times \mathbb{S}^{1}) > 0 \Rightarrow \lambda_{\tau^{*}h}(X) > 0.
\end{equation*}
\end{corollary}

\begin{proof}  If $\lambda_{h}(X \times \mathbb{S}^{1}) > 0$, there is a metric $ h'$ in the conformal class $[h]$ of constant positive scalar curvature, by the solution of the Yamabe problem.  The condition (\ref{eq:intro}) is conformally invariant, so $h'$ also satisfies (\ref{eq:intro}).  Since the
proof of Thm.~ \ref{MAIN:thm1} involves the conformal transformation $g = h' + dt^2\mapsto \tilde g$, $X$ admits a PSC metric in $[\tau^*\imath^*g]
= [\tau^* h'] =[\tau^* h].$  Since the sign of the Yamabe invariant is conformally invariant, $\lambda_{\tau^*h}(X)>0.$
\end{proof}

\bibliographystyle{plain}
\bibliography{S1Stability}

\begin{thebibliography}{10}

\bibitem{Niren4}
S.~Agmon, A.~Douglis, and L.~Nirenberg.
\newblock Estimates near the boundary for solutions of elliptic partial
  differential equstions satisfying general boundary conditions {I}.
\newblock {\em Commun. Pure Appl. Math}, 12:623--727, 1959.

\bibitem{AB2}
K.~Akutagawa and B.~Botvinnik.
\newblock The relative {Y}amabe invariant.
\newblock {\em Comm. Anal. Geom.}, 10(5):935--965, 2002.

\bibitem{carlotto}
A.~Carlotto and C.~Li.
\newblock Constrained deformations of positive scalar curvature metrics.
\newblock {\em J. Differential Geometry}, 126(2), 2024.

\bibitem{Che}
P.~Cherrier.
\newblock Probl\'emes de {N}eumann non lin\'eaires sur les vari\'et\'es
  {R}iemanniannes.
\newblock {\em J. Funct. Anal.}, 57:154--206, 1984.

\bibitem{Chodosh}
O.~Chodosh.
\newblock Stable minimal surfaces and positive scalar curvature.
  \url{https://web.stanford.edu/~ochodosh/Math258-min-surf.pdf}.

\bibitem{KS}
A.~Kumar and B.~Sen.
\newblock Positive scalar curvature and exotic structures on simply connected
  four manifolds. \url{https://arxiv.org/abs/2501.01113}.
\newblock 2025.

\bibitem{Rade}
D.~R\"ade.
\newblock Scalar and mean curvature comparison via $ \mu $-bubbles.
\newblock {\em Calc. Var. Partial Differential Equations}, 62(187), 2023.

\bibitem{JR}
J.~Rosenberg.
\newblock Manifolds of positive scalar curvature: a progress report.
\newblock {\em Surveys in Differential Geometry}, 11(1):259--294, 2006.

\bibitem{RRX}
D.~Ruberman, S.~Rosenberg, and J.~Xu.
\newblock The conformal {L}aplacian and positive scalar curvature metrics on
  manifolds with boundary. {\url{https://https://arxiv.org/abs/2302.05521}}.
\newblock 2023.

\bibitem{Zeidler}
R.~Zeidler.
\newblock Band width estimates via the {D}irac operator.
\newblock {\em J. Differential Geom.}, 122(1):155--183, 2022.

\end{thebibliography}

\end{document}